\documentclass[11pt,a4paper]{article}

\usepackage{caption2}
\usepackage{CJK}
\usepackage{tabls}
\usepackage{bm}
\usepackage{multirow}
\usepackage{amssymb}
\usepackage{amsfonts}
\usepackage{mathrsfs}
\usepackage{empheq}
\usepackage{amsmath}
\usepackage{amsthm}
\usepackage{graphicx}
\usepackage{color}
\usepackage{indentfirst}
\usepackage{hyperref}
\usepackage{float}
\usepackage{cases}
\usepackage[titletoc]{appendix}
\usepackage{bm}

\usepackage[colorinlistoftodos,english]{todonotes}

\allowdisplaybreaks%

\def\i1n{i=1,\cdots,n}
\def\j1n{j=1,\cdots,n}
\def\ij1n{i,j=1,\cdots,n}

\def \i{\mathrm i}

 \numberwithin{equation}{section}
\theoremstyle{definition}

 \newtheorem{thm}{\indent Theorem}[section]
 
 \newtheorem{lem}{\indent Lemma}[section]

 \newtheorem{theorem}{Theorem}[section]

 \newtheorem{remark}[theorem]{Remark}

\theoremstyle{definition}
\theoremstyle{theorem}
\theoremstyle{lemma}

\newcommand{\be}{\begin{equation}}
\newcommand{\ee}{\end{equation}}
\newcommand{\beq}{\begin{equation*}}
\newcommand{\eeq}{\end{equation*}}

\begin{document}

\begin{CJK*}{GB}{gbsn}
\title{1+2 dimensional radially symmetric wave maps revisit.}

\author{Yi Zhou \thanks{School of Mathematics Science, Fudan University, Shanghai, P. R. China (yizhou@fudan.edu.cn).}}

\date{}

\pagestyle{myheadings} \markboth{1+2 dimensional radially symmetric wave maps revisit.}{1+2 dimensional radially symmetric wave maps revisit.}

\maketitle

\begin{abstract}
The author gives an alternative and simple proof of the global existence of smooth solutions to the Cauchy problem for wave maps from the 1+2-dimensional Minkowski space to an arbitrary compact smooth Riemannian manifold without boundary, for arbitrary smooth, radially symmetric data. the author can also treat non-compact manifold under some additional assumptions which generalize the existing ones.
\end{abstract}


\section{Main result} \label{section1}

Let $N$ be a smooth Riemannian k-manifold without boundary. With no loss of generality, we assume that $ N \subset \mathbb{R}^{n}$, isometrically. We consider wave maps $ \Phi = ( \Phi^{1}, \cdots, \Phi^{n} )=\Phi(t,x): \mathbb{R} \times \mathbb{R}^{2} \rightarrow N \in \mathbb{R}^{n}  $ , satisfying the equation
\begin{equation}
 \Box \Phi = \Phi_{tt} - \Delta \Phi = B(\Phi)( \partial_{\alpha} \Phi, \partial^{\alpha} \Phi ) \perp T_{\Phi} N ,
\end{equation}
where $B$ denotes the second fundamental form of $N$. Writing $z = (t,x) = (x^{\alpha})_{0 \leq \alpha \leq 2}$, we also let $ \partial_{\alpha}= \frac{\partial}{\partial x^{\alpha}}, \ \alpha=0,1,2. $ We raise and lower indices with the Minkowski metric $ \eta = (\eta_{\alpha \beta}) = ( \eta^{\alpha \beta} ) = diag(-1,1,1) $ and tacitly sum over repeated indices.

Due to its mathematical difficulty and important physical background, the topic of wave maps has experienced an incredible advancement in the past several decades. It has at least two physical motivations to study wave maps. One is the nonlinear $\sigma $-model which deals with the case $N$ is a sphere and on the other hand, vacuum Einstein equations with $U(1) \times \mathbb{R} $ symmetries reduce to a radially symmetric wave maps from 1+2 dimensional Minkowski space to the hyperbolic plane. It is Prof. Gu \cite{3} who first gave the regularity result in $1+1$ dimensional case. For an up to date account of the full developments,  we refer to the monograph of Geba, D.A. and Grillakis M.G. \cite{2}, as well as
\cite{4}, \cite{5}, and \cite{8}.

The main purpose of this paper is to give an alternative proof of the global existence for 1+2 dimensional wave maps with radial symmetry.
We may assume that TN is parallelizable. Let $ \bar{e}_{1} , \cdots , \bar{e}_{k} $ be a smooth orthonormal frame field such that at any point $ p \in N $ the vectors $ \bar{e}_{1}(p) , \cdots , \bar{e}_{k}(p) $ form an orthonormal basis for $ T_{p}N $.We make the following assumption on the target manifold.

(H1): $ \nabla_{p} \bar{e}_{i}(p) $ is uniformly bounded,

(H2): $ \nabla_{p} B(p) $ is uniformly bounded.

We note that those assumptions are in particular satisfied for compact manifold.

The main result of this paper is the following.

\begin{thm}
Let $ N \subset \mathbb{R}^{n}$ be a smooth Riemannian manifold without boundary. Assume (H1) (H2) satisfied for $N$. Then for any radially symmetric data
\begin{equation}
 ( \Phi_{0}, \Phi_{1} ) = ( \Phi_{0}(r), \Phi_{1}(r) ) \in C^{\infty} (\mathbb{R}^{2}, TN), \ r=|x|,
\end{equation}	
there exists a unique, smooth solution $\Phi=\Phi(t,r) $ to the Cauchy problem (1.1)(1.2), defined for all time.
\end{thm}

Our result slightly generalizes the work of Christodoulou-Tahvildar-Zadeh \cite{1} and Struwe \cite{6},\cite{7}.

As usual, the proof of Theorem1.1 is divided  into two main steps. The first step is to show that small energy implies regularity and the second step is to show that energy can not concentrate. In this paper, we give an alternative proof of the first step which is totally different from the work of Christodoulou-Tahvildar-Zadeh \cite{1} and the proof of the second step can be found in Struwe \cite{6},\cite{7}. In the work of Christodoulou-Tahvildar-Zadeh \cite{1}, the first step is achieved by a H$\ddot{o}$lder estimate using the fundamental solution of the 2-dimensional radially symmetric wave operator which is quite complicated, While we rely on an energetic argument which we call a new div curl lemma  and has potential to work for quasilinear problems. This will be pursued in our future work.

\section{Intrinsic setting}

Let $\Phi = \Phi(t,r) $ be a smooth radially symmetric wave maps, then (1.1)(1.2) can be rewriting as
\begin{align}
\Phi_{tt} - \Phi_{rr} - \frac{\Phi_{r}}{r} = 4 B(\Phi)(\Phi_{u}, \Phi_{v}), \ (t>0)
\end{align}

\begin{align}
\Phi= \Phi_{0}(r), \ \Phi_{t}= \Phi_{1}(r), \ (t=0)
\end{align}
we denote  $ u=t-r, \ v=t+r, $ then
\begin{align*}
 &t=\frac{1}{2}(u+v),
 \partial_{t}=\partial_{u}+\partial_{v},  \partial_{u}=\frac{1}{2}(\partial_{t}-\partial_{r}), \\
 &r=\frac{1}{2}(v-u), \partial_{r}=\partial_{v}-\partial_{u}, \partial_{v}=\frac{1}{2}(\partial_{t}+\partial_{r}), \\
 &\Phi_{u}=\partial_{u} \Phi, \Phi_{v}=\partial_{v} \Phi, \Phi_{tt}=\partial_{t}^{2} \Phi, \ etc.
\end{align*}

and we have the energy conservation law
\begin{align}
E(t) = \frac{1}{2} \int_{0}^{\infty} r(|\Phi_{t}|^{2} + |\Phi_{r}|^{2} ) dr = E_{0} = \frac{1}{2} \int_{0}^{\infty} r(|\Phi_{0}^{\prime}|^{2} + |\Phi_{1}|^{2} ) dr.
\end{align}

Let $D$ be the pull-back covariant derivative in $u^{*} TN $, we may write equation(2.1) as
\begin{align}
D_{t} \Phi_{t}-\frac{1}{r}D_{r}(r \Phi_{r}) = 0.
\end{align}
 From $ (\bar{e}_{i} )_{1 \leq i \leq k } $, we obtain a frame $ e_{i} = R_{ij}( \bar{e}_{j} \circ \Phi ), 1 \leq i \leq k $ for the pull-back bundle, where $R=R(z) = (R_{ij}) $ may be any smooth map. Denoting
\begin{align}
D_{\alpha}e_{i}=A_{i \alpha }^{j} e_{j}, \ 0 \leq \alpha \leq 2
\end{align}
with a matrix-valued connection 1-form $ A=A_{\alpha} dx^{\alpha} $. We compute the curvature $F$ of $D$ via the commutation relation
\begin{align*}
& D_{\alpha}D_{\beta}e_{i} - D_{\beta}D_{\alpha}e_{i} = D_{\alpha}(A_{i \beta }^{j} e_{j}) -D_{\beta}(A_{i \alpha }^{j} e_{i}) \\
& = (\partial_{\alpha}A_{i \beta }^{k} - \partial_{\beta}A_{i \alpha }^{k} - A_{j \alpha }^{k} A_{i \beta }^{j} - A_{j \beta }^{k} A_{i \alpha }^{j} ) e_{k} \\
& = F_{i \alpha \beta }^{k} e_{k}
\end{align*}
or, more precisely
\begin{align}
dA + \frac{1}{2} [A,A]=F.
\end{align}
In the radially symmetric case, we may choose $ R=R(t,r), A = A_{0}dt+A_{1}dr $. Following Struwe \cite{7},  we may impose "exponential gauge" condition $A_{1}=0$, which yields the relation
\begin{align}
*dA = -\partial_{r}A_{0} = F_{01}.
\end{align}
If we normalize $ A_{0}(t, \infty )=0, \forall t $, from this relation we obtain
\begin{align}
A_{0} = \int_{r}^{+\infty} F_{01} ds.
\end{align}
By (H1), we get
\begin{align}
|F_{01}| \leq C| d \Phi|^{2}.
\end{align}
Thus we deduce the estimate
\begin{align}
|A_{0}| \lesssim \int_{r}^{+\infty} | d \Phi|^{2} ds \lesssim E_{0}r^{-1}.
\end{align}
Let
\begin{align}
\Phi_{t} = q_{0}^{i}e_{i}, \ \Phi_{r} = q_{1}^{i}e_{i}.
\end{align}
Using the notation
\begin{align}
D_{\alpha}\partial_{\beta} \Phi &=D_{\alpha}(q_{\beta}^{i}e_{i}) = (\partial_{\alpha}q_{\beta}^{j} + A_{i \beta}^{j}q_{\beta}^{i} )e_{j} \notag \\
&=(D_{\alpha}q_{\beta}^{i})^{j}e_{j},
\end{align}
we then may write equation (2.4) in the form
\begin{align}\label{1}
D_{t}q_{0} - \frac{1}{r}D_{r}(rq_{1}) = \partial_{t} q_{0} + A_{0}q_{0} - \frac{1}{r}\partial_{r}(rq_{1}) =0.
\end{align}
Moreover, we have the commutation relation
\begin{align}\label{2}
D_{r}q_{0}= D_{t}q_{1} = \partial_{t} q_{1} + A_{0}q_{1}.
\end{align}

Taking inner product in $ \mathbb{R}^{n} $ of \eqref{1} with $-r^{1-\alpha}q_{1}$, we get
\begin{align} \label{2022_2_3}
0&=-r^{1-\alpha} q_1 \partial_{t} q_0-r^{1-\alpha} q_1 (A_0 q_0)+r^{1-\alpha} q_1\partial_r q_1 +r^{-\alpha} q_1^2\nonumber\\
&=-r^{1-\alpha} q_1 \partial_{t} q_0-r^{1-\alpha} q_1 (A_0 q_0)+ \left( r^{1-\alpha} q_1\partial_r q_1+\frac12 (1-\alpha)r^{-\alpha} q_1^2\right)+\frac12 (1+\alpha)r^{-\alpha} q_1^2\nonumber\\
&=-r^{1-\alpha} q_1 \partial_{t} q_0-r^{1-\alpha} q_1 (A_0 q_0)+ \partial_r\left( \frac12 r^{1-\alpha} q_1^2\right)+\frac12 (1+\alpha)r^{-\alpha} q_1^2.
\end{align}

Taking inner product in $ \mathbb{R}^{n} $ of \eqref{2} with $-r^{1-\alpha}q_{0}$, we get
\begin{align}\label{2022_2_4}
0&=	-r^{1-\alpha} q_0\partial_{t} q_1	-r^{1-\alpha} q_0(A_0 q_1)+	r^{1-\alpha} q_0 \partial_r q_0 \nonumber\\
&=	-r^{1-\alpha} q_0\partial_{t} q_1	-r^{1-\alpha} q_0(A_0 q_1)+\left( r^{1-\alpha} q_0\partial_r q_0+\frac12 (1-\alpha)r^{-\alpha} q_0^2\right)- \frac12 (1-\alpha)r^{-\alpha} q_0^2\nonumber\\
&=	-r^{1-\alpha} q_0\partial_{t} q_1	-r^{1-\alpha} q_0(A_0 q_1)+\partial_r\left(\frac12 r^{1-\alpha} q_0^2\right)- \frac12 (1-\alpha)r^{-\alpha} q_0^2.
\end{align}

Combining \eqref{2022_2_3}  \eqref{2022_2_4} and the antisymmetry of $A_0$,  we know
\begin{align} \label{2022_2_5}
&	-\partial_{t} \left( r^{1-\alpha} q_0 q_1 \right)+\partial_r\left(\frac12 r^{1-\alpha} (q_0^2+q_1^2)\right)- \frac12 (1-\alpha)r^{-\alpha} q_0^2+\frac12 (1+\alpha)r^{-\alpha} q_1^2\nonumber\\
&	=r^{1-\alpha} q_1 (A_0 q_0)+r^{1-\alpha} q_0(A_0 q_1) = 0 .
\end{align}

Define $ Q_0(r)= \int_0^r \xi^{-\sigma} q_1(\xi) d \xi $, $ Q(r)=r^{-\alpha+\sigma} Q_0 $,   $ \sigma<\alpha $.  And  by \eqref{2}, we have
\begin{align} \label{6}
Q \partial_{t} q_0 &=\partial_{t} (Q q_0) - \partial_{t} Q q_0=\partial_{t} (r^{-\alpha+\sigma} Q_0 q_0) -  r^{-\alpha+\sigma} \partial_{t} Q_0 q_0\nonumber\\
&=\partial_{t} (r^{-\alpha+\sigma} Q_0 q_0) - \frac{q_0}{r^{\alpha-\sigma}}  \int_0^r \xi^{-\sigma} \partial_{t} q_1(\xi) d \xi \nonumber \\
&= \partial_{t} (r^{-\alpha+\sigma} Q_0 q_0) -  \frac{q_0}{r^{\alpha-\sigma}}\int_0^r \xi^{-\sigma} [\partial_{\xi}  q_0(\xi)- A_0 q_1] d \xi\nonumber\\
&= \partial_{t} (r^{-\alpha+\sigma} Q_0 q_0) +\frac{q_0}{r^{\alpha-\sigma}}\int_0^r \xi^{-\sigma} A_0 q_1 d \xi-   \frac{q_0(r)}{r^{\alpha-\sigma}}\left[ \left.\xi^{-\sigma}  q_0(\xi)\right|_{0}^r+\sigma \int _0^r \xi^{-\sigma-1} q_0(\xi) d \xi\right] \nonumber\\
&= \partial_{t} (r^{-\alpha+\sigma} Q_0 q_0) +\frac{q_0}{r^{\alpha-\sigma}}\int_0^r \xi^{-\sigma} A_0 q_1 d \xi-    \frac {q_0^2}{r^\alpha}- \frac{\sigma q_0}{r^{\alpha-\sigma}} \int _0^r \xi^{-\sigma-1} q_0(\xi) d \xi
\end{align}

\begin{align} \label{7}
-\partial_r q_1 Q &=-\partial_r (q_1 Q)+q_1 \partial_r Q=-\partial_r (q_1 Q)+q_1 \partial_r (r^{-\alpha+\sigma} Q_0)\nonumber\\
&= -\partial_r (q_1 Q)+ (-\alpha+\sigma ) r^{-\alpha+\sigma-1} q_1  Q_0+r^{-\alpha+\sigma} q_1 \partial_r  Q_0\nonumber\\
&=-\partial_r (q_1 r^{-\alpha+\sigma} Q_0)+ (-\alpha+\sigma ) \frac {q_1  Q}{r}+\dfrac {q_1^2} {r^\alpha}
\end{align}
\begin{align} \label{8}
\frac{q_1}{r} Q&=  \frac{q_1}{r^{\alpha+1-\sigma}} Q_0= \frac{1}{r^{\alpha+1-2\sigma}} \partial_r  Q_0  Q_0\nonumber\\
&= \frac12 \partial_r (r^{2\sigma-1-\alpha}   Q_0 ^2)-\frac{2\sigma-1-\alpha}{2}r^{2\sigma-2-\alpha} Q_0^2
\end{align}

Combining \eqref{7}£¬\eqref{8}, we can get
\begin{align} \label{9}
-\partial_r q_1 Q	-\frac{q_1}{r} Q&=-\partial_r (q_1 r^{-\alpha+\sigma} Q_0)+ (-\alpha+\sigma -1) \frac {q_1  Q}{r}+\dfrac {q_1^2} {r^\alpha}\nonumber\\
&=-\partial_r (q_1 r^{-\alpha+\sigma} Q_0)+\dfrac {q_1^2} {r^\alpha}\nonumber\\
&+ \frac{-\alpha+\sigma -1}{2} \partial_r (r^{2\sigma-1-\alpha}   Q_0 ^2)-\frac{(2\sigma-1-\alpha)(-\alpha+\sigma -1)}{2}r^{2\sigma-2-\alpha} Q_0^2
\end{align}

From \eqref{6} and \eqref{9}, taking inner product in $ \mathbb{R}^{n} $ of \eqref{1} with	$ -\dfrac{Q}{2}$, we know
\begin{align}  \label{10}
-\frac 12 &\partial_{t} (r^{-\alpha+\sigma} Q_0 q_0) +\frac 12  \frac {q_0^2}{r^\alpha}-\frac 12 \dfrac {q_1^2} {r^\alpha}+\frac 12 \partial_r (q_1 r^{-\alpha+\sigma} Q_0)\nonumber\\
&+\frac{\alpha-\sigma +1}{4} \partial_r (r^{2\sigma-1-\alpha}   Q_0 ^2)+\frac{(1+\alpha-2\sigma)(\alpha-\sigma +1)}{4}r^{2\sigma-2-\alpha} Q_0^2 \nonumber\\
&= \frac{q_0}{2r^{\alpha-\sigma}}\int_0^r \xi^{-\sigma} A_0 q_1 d \xi- \frac{\sigma q_0}{2r^{\alpha-\sigma}} \int _0^r \xi^{-\sigma-1} q_0(\xi) d \xi+ \frac 12  r^{-\alpha+\sigma} A_0q_0  Q_0
\end{align}

By  \eqref{2022_2_5}  and \eqref{10} ,  taking inner product in $ \mathbb{R}^{n} $ of \eqref{1} with
$ -(\dfrac{Q}{2} +r^{1-\alpha} q_1) $,  we obtain
\begin{align} \label{11}
&	-\partial_{t} \left( r^{1-\alpha} q_0 q_1 \right)+\partial_r\left(\frac12 r^{1-\alpha} (q_0^2+q_1^2)\right)+\frac\alpha2 \frac{ q_0^2+q_1^2}{r^{\alpha}} \nonumber\\
&		-\frac 12 \partial_{t} (r^{-\alpha+\sigma} Q_0 q_0) +\frac 12 \partial_r (q_1 r^{-\alpha+\sigma} Q_0)\nonumber\\
&+\frac{\alpha-\sigma +1}{4} \partial_r (r^{2\sigma-1-\alpha}   Q_0 ^2)+ \frac{(1+\alpha-2\sigma)(\alpha-\sigma +1)}{4}r^{2\sigma-2-\alpha} Q_0^2 \nonumber\\
&	= \frac{q_0}{2r^{\alpha-\sigma}}\int_0^r \xi^{-\sigma} A_0 q_1 d \xi- \frac{\sigma q_0}{2r^{\alpha-\sigma}} \int _0^r \xi^{-\sigma-1} q_0(\xi) d \xi+ \frac 12  r^{-\alpha+\sigma} A_0q_0  Q_0.
\end{align}

\section{New div curl Lemma}

The purpose of this section is to prove the following lemma3.1 and lemma3.2. We call Lemma3.2 a new div curl Lemma.

\begin{lem}
Suppose that
\begin{equation}
  \begin{cases}
  \partial_{u} F^{11} + \partial_{v} F^{12} = G^{1} \ \ (r\geq 0, 0 \leq t \leq T) \\
  \partial_{u} F^{21} - \partial_{v} F^{22} = G^{2} \ \ (r\geq 0, 0 \leq t \leq T),
  \end{cases}
\end{equation}
where $ F^{11},F^{12},F^{21},F^{22} $ are all nonnegative. Moreover,
	 \begin{align}
	r=0: \ F^{11}-F^{12}=F^{21}+F^{22}=0.
	 \end{align}
Then there hold
\begin{align}
& \sup_{u} \int_{|u|}^{2T-u} F^{11}(u,v) du + \sup_{v} \int_{-v}^{min\{v,2T-v\}} F^{12}(u,v) du \notag \\
& \lesssim \int_{0}^{\infty} (F^{11}+F^{12}) (0,r) dr +
\int_{0}^{T} \int_{0}^{\infty} |G_{1}| (t,r) drdt,\\
& \sup_{u} \int_{|u|}^{2T-u} F^{21}(u,v) du + \sup_{v} \int_{-v}^{min\{v,2T-v\}} F^{22}(u,v) du \notag \\
& \lesssim \int_{0}^{\infty} (F^{21}+F^{22}) (T,r) dr + \int_{0}^{\infty} (F^{21}+F^{22}) (0,r) dr +
\int_{0}^{T} \int_{0}^{\infty} |G_{2}| (t,r) drdt. \notag \\
\end{align}	
\end{lem}

{\bf Proof}\ \ We only prove
\begin{align*}
& \sup_{u} \int_{|u|}^{2T-u} F^{11}(u,v) du  \\
& \lesssim \int_{0}^{\infty} (F^{11}+F^{12}) (0,r) dr +
\int_{0}^{T} \int_{0}^{\infty} |G_{1}| (t,r) drdt.
\end{align*}	
The other estimates are similar. Noting(3.2), we have
\begin{align*}
& \partial_{u} \int_{|u|}^{2T-u} F^{11}(u,v) du  \\
& = -F^{11}|_{v=2T-u} - F^{11}|_{v=|u|} sgn(u) +  \int_{|u|}^{2T-u} \partial_{u} F^{11} dv \\
& = -F^{11}|_{v=2T-u} - F^{11}|_{v=|u|} sgn(u) - \int_{|u|}^{2T-u} \partial_{v} F^{12} dv +
\int_{|u|}^{2T-u} |G_{1}| dv \\
& = -F^{11}|_{v=2T-u} -F^{12}|_{v=2T-u} - F^{11}|_{v=|u|} sgn(u) +F^{12}|_{v=|u|} +
\int_{|u|}^{2T-u} |G_{1}| dv \\
& \lesssim (F^{11}+F^{12})|_{t=0} +
\int_{|u|}^{2T-u} |G_{1}| dv.
\end{align*}	
Integration in $u$ yields the desired estimate.

\begin{lem}
	
	Under the assumption of Lemma3.1, we have
\begin{align*}
& \int_{0}^{T} \int_{0}^{\infty} (F^{11}F^{22}+F^{12}F^{21}) (t,r) drdt \\
& \lesssim \left( \int_{0}^{\infty}(F^{11}+F^{12}) (0,r) dr + \int_{0}^{T} \int_{0}^{\infty} |G_{1}(t,r)| drdt \right)  \\
& \cdot \left( \int_{0}^{\infty}(F^{21}+F^{22}) (0,r) dr + \int_{0}^{\infty}(F^{21}+F^{22}) (T,r) dr + \int_{0}^{T} \int_{0}^{\infty} |G_{2}(t,r)| drdt \right) .
\end{align*}	
\end{lem}

{\bf Proof}\ \ Similary to the proof of Lemma3.1, for $ \bar{v} \leq 2T-u $, we get
\begin{align*}
&  \partial_{u} \int_{|u|}^{\bar{v}} F^{11}(u,v) dv + F^{12}(u,\bar{v}) \\
& \lesssim (F^{11}+F^{12})|_{u+\bar{v}=0} + \int_{|u|}^{\bar{v}} G_{1} dv.
\end{align*}
Thus
\begin{align*}
& F^{21}(u,v) \partial_{u} \int_{|u|}^{\bar{v}} F^{11}(u,v) dv + F^{21}(u,\bar{v}) F^{12}(u,\bar{v})  \\
& \lesssim F^{21}(u,\bar{v}) (F^{11}+F^{12})|_{u+\bar{v}=0} + \left( \int_{|u|}^{\bar{v}} G_{1}(u,v) dv \right) F^{21}(u,\bar{v}).
\end{align*}
By the second equation in (3.1), the first term is
\begin{align*}
& \partial_{u} \left[ F^{21}(u,\bar{v}) \int_{|u|}^{\bar{v}} F^{11}(u,v) dv  \right]
- \partial_{u}  F^{21}(u,\bar{v}) \int_{|u|}^{\bar{v}} F^{11}(u,v) dv  \\
& = \partial_{u} \left[ F^{21}(u,\bar{v}) \int_{|u|}^{\bar{v}} F^{11}(u,v) dv  \right]
-\partial_{\bar{v}} F^{22}(u,\bar{v}) \int_{|u|}^{\bar{v}} F^{11}(u,v) dv \\
& - G_{2}(u,\bar{v}) \int_{|u|}^{\bar{v}} F^{11}(u,v) dv \\
& = \partial_{u} \left[ F^{21}(u,\bar{v}) \int_{|u|}^{\bar{v}} F^{11}(u,v) dv  \right]
- \partial_{\bar{v}} \left[ F^{22}(u,\bar{v}) \int_{|u|}^{\bar{v}} F^{11}(u,v) dv  \right]  \\
& + F^{22}(u,\bar{v}) F^{11}(u,\bar{v})
- G_{2}(u,\bar{v}) \int_{|u|}^{\bar{v}} F^{11}(u,v) dv.
\end{align*}
Therefore, we get
\begin{align*}
& F^{22}(u,\bar{v}) F^{11}(u,\bar{v}) + F^{21}(u,\bar{v}) F^{12}(u,\bar{v}) \\
& \lesssim \partial_{\bar{v}} \left[ F^{22}(u,\bar{v}) \int_{|u|}^{\bar{v}} F^{11}(u,v) dv  \right]
- \partial_{u} \left[ F^{21}(u,\bar{v}) \int_{|u|}^{\bar{v}} F^{11}(u,v) dv  \right] \\
& + G_{2}(u,\bar{v}) \int_{|u|}^{\bar{v}} F^{11}(u,v) dv
+ F^{21}(u,\bar{v})(F^{11}+F^{12})|_{u+\bar{v}=0}
+ \left( \int_{|u|}^{\bar{v}} G_{1}(u,v) dv  \right) F^{21}(u,\bar{v}).
\end{align*}
Integration in the region $ \left\lbrace 0 \leq u+\bar{v} \leq 2T, \ \bar{v}-u \geq 0 \right\rbrace  $ and using Lemma3.1 for the boundary estimate, we get the desired conclusions.

\section{Small energy implies regularity}

In this section, we assume
\begin{align}
E_{0} \leq \epsilon_{1}
\end{align}
is sufficiently small. Under this assumption, we shall prove the solution is smooth. For that purpose, we only need to give the $H^{2} $ estimate of the solution, see \cite{4}.

Let $ \Psi = \Phi_{t} $. Differentiating the equation(2.1), we get
\begin{align}
& \Psi_{tt} - \Psi_{rr} - \frac{\Psi_{r}}{r} \notag \\
& = 4 B(\Phi) ( \Psi_{u}, \Phi_{v} )
+ 4 B(\Phi) (\Phi_{u}, \Psi_{v}  )
+ 4 \Psi \cdot B^{\prime}(\Phi) ( \Phi_{u}, \Phi_{v} ),
\end{align}
thus
\begin{align*}
& r\Psi_{t} \cdot \left( \Psi_{tt} - \Psi_{rr} - \frac{\Psi_{r}}{r} \right)  \\
& = 4r\Psi_{t} \cdot B(\Phi) ( \Psi_{u}, \Phi_{v} )
+ 4r\Psi_{t} \cdot B(\Phi) (\Phi_{u}, \Psi_{v}  )
+ 4\Psi_{t} \cdot \left[  \Psi \cdot A^{\prime}(\Phi) ( \Phi_{u}, \Phi_{v} ) \right] , \\
& \Psi \cdot B(\Phi) =0.
\end{align*}
Differentiating with respect to $t$ yields
\begin{align*}
\Psi_{t} \cdot B(\Phi) = - \Psi \cdot \left[ \Psi \cdot A^{\prime}(\Phi) \right].
\end{align*}
So we get
\begin{align}
& \partial_{t} \left[ \frac{1}{2} (\Psi_{t}^{2}+ \Psi_{r}^{2} ) r  \right]
- \partial_{r} \left[ r \Psi_{t} \cdot  \Psi_{r} \right] \notag  \\
& = 4r\Psi \cdot \left[ \Psi \cdot B^{\prime}(\Phi) \right]  ( \Psi_{u}, \Phi_{v} )
+ 4r\Psi \cdot \left[ \Psi \cdot B^{\prime}(\Phi) \right] (\Phi_{u}, \Psi_{v}  )
+ 4r\Psi_{t} \cdot \left[ \Psi \cdot B^{\prime}(\Phi) \right] ( \Phi_{u}, \Phi_{v} ) \notag \\
& =G_{1}.
\end{align}
In a similar way, let $ \hat{\Psi} = \Phi_{r}(t,r) $. We get
\begin{align}
& \partial_{t} \left[ \frac{1}{2} (\hat{\Psi}_{t}^{2}+ \hat{\Psi}_{r}^{2} + \frac{\hat{\Psi}^{2}}{r^{2}} ) r  \right]
- \partial_{r} \left[ r \hat{\Psi}_{t} \cdot  \hat{\Psi}_{r} \right] \notag  \\
& = 4r \hat{\Psi} \cdot \left[ \hat{\Psi} \cdot B^{\prime}(\Phi) \right]  ( \hat{\Psi}_{u}, \Phi_{v} )
+ 4r\hat{\Psi} \cdot \left[ \hat{\Psi} \cdot B^{\prime}(\Phi) \right] (\Phi_{u}, \hat{\Psi}_{v}  )
+ 4r\hat{\Psi}_{t} \cdot \left[ \hat{\Psi} \cdot B^{\prime}(\Phi) \right] ( \Phi_{u}, \Phi_{v} ) \notag \\
& =\tilde{G}_{1}.
\end{align}
Integrating in $t,r$, we get
\begin{align}
 \sup_{0 \leq t \leq T} | D^{2} \Phi |^{2}_{L^{2}(\mathbb{R}^{2})} (t)
 \lesssim | \nabla^{2} u_{0} |^{2}_{L^{2}(\mathbb{R}^{2})}
 +| \nabla u_{1} |^{2}_{L^{2}(\mathbb{R}^{2})}
 + \int_{0}^{T}  \int_{0}^{\infty} ( |\tilde{G}_{1}| + |G_{1}| ) (t,r) drdt.
\end{align}

The main purpose of this section is to prove
\begin{align}
\int_{0}^{T}  \int_{0}^{\infty} ( |\tilde{G}_{1}| + |G_{1}| ) (t,r) drdt
\lesssim \epsilon_{1} E_{1} ,
\end{align}
where
\begin{align}
E_{1} =  | \nabla^{2} u_{0} |^{2}_{L^{2}(\mathbb{R}^{2})}
+| \nabla u_{1} |^{2}_{L^{2}(\mathbb{R}^{2})}.
\end{align}
For that purpose, we use an induction argument and first we assume
\begin{align}
\int_{0}^{T}  \int_{0}^{\infty} ( |\tilde{G}_{1}| + |G_{1}| ) (t,r) drdt
\lesssim E_{1}.
\end{align}
Then we get
\begin{align}
| D^{2} \Phi |^{2}_{L^{2}(\mathbb{R}^{2})} (t)
\lesssim E_{1}.
\end{align}
To obtain (4.6), we only estimate $G_{1}$. The estimate of $ \tilde{G}_{1}  $ can be done in a similar way. By the expression of $G_{1}$, we get
\begin{align}
|G_{1}| & \lesssim | \Psi |^{2} ( | \Psi_{u} | | \Phi _{v}| + | \Phi_{u}| | \Psi_{v} | )r
+ | \Psi_{t} | ( | \Psi | | \Phi_{u}| | \Phi_{v}| )r \notag \\
& = g_{1}+g_{2}.
\end{align}
We first estimate $g_{1}$. Let $ \frac{1}{10} \leq \beta \leq \frac{1}{4} $, we have
\begin{align}
\int_{0}^{T}  \int_{0}^{\infty} |g_{1}| drdt
 \lesssim \left[   \int_{0}^{T}  \int_{0}^{\infty} \left( | \Psi_{u} |^{2} | \Phi_{v}|^{2} + | \Phi_{u}|^{2} | \Psi_{v} |^{2}  \right) r^{2-\beta} drdt \right]^{\frac{1}{2} }
  \left[   \int_{0}^{T}  \int_{0}^{\infty} | \Psi |^{4}  r^{\beta} drdt \right]^{\frac{1}{2} }.
\end{align}
By Sobolev-Hardy inequality,
\begin{align}
\int_{0}^{\infty} | \Psi |^{4}  r^{\beta} dr \lesssim \left( \int_{0}^{\infty} | \Psi |^{2}  r^{-\beta} dr \right)
\left( \int_{0}^{\infty} | \Psi |^{2}  r dr \right)^{\beta}
 \left( \int_{0}^{\infty} | \Psi_{r} |^{2}  r dr \right)^{1-\beta},
\end{align}
thus
\begin{align}
\int_{0}^{T}  \int_{0}^{\infty} |g_{1}| drdt
& \lesssim \left[   \int_{0}^{T}  \int_{0}^{\infty} \left( | \Psi_{u} |^{2} | \Phi_{v}|^{2} + | \Phi_{u}|^{2} | \Psi_{v} |^{2}  \right) r^{2-\beta} drdt \right]^{\frac{1}{2} }
\left[   \int_{0}^{T}  \int_{0}^{\infty} | \Psi |^{2}  r^{-\beta} drdt \right]^{\frac{1}{2} } \notag  \\
\end{align}
Take $ \alpha = \beta $ in \eqref{11} and integrate for $ 0 \leq r \leq +\infty, \  0 \leq t \leq T $. Noting (2.19), we get
\begin{align*}
& \frac{1}{2} \int_{0}^{T}  \int_{0}^{\infty} r^{-\beta} \left( | \Psi |^{2} + |\hat{\Psi}|^{2}  \right) drdt
= \left( \int_{0}^{\infty}  r^{1-\beta} \Psi \hat{\Psi} \right)|_{0}^{T}
\left(q^2_0+q^2_1\right)drdt \notag \\
& -\frac{1}{2}\int_{0}^{\infty} r^{-\beta+\sigma} Q_0 q_0 dr|_{0}^{T}\notag \\
& +\left[\frac{1}{2}\int_{0}^{\infty}q_0r^{\sigma -\alpha}\int_{0}^{r} \xi^{\alpha}A_0q_1d\xi-\sigma q_0 r^{\sigma-\alpha} \int_{0}^{r} \xi^{-\sigma-1}q_0(\xi)d\xi +\frac{1}{2} r^{\sigma-\alpha}A_0q_0Q_0 drdt\right]
\end{align*}

By Hardy's inequality and noting (2.10), the third term in the right hand of side of the above equality can be absorbed by the left hand side, and by Sobolev-hardy inequality the first and second term is bounded by
\begin{align*}
	\epsilon_{1}^{1-\beta/2}E^{\beta/2}_1,
\end{align*}
thus
\begin{align}
\left(  \int_{0}^{T}  \int_{0}^{\infty} r^{-\beta} \left( | \Psi |^{2} + |\hat{\Psi}|^{2}  \right) drdt \right)^{\frac{1}{2}}
\lesssim \epsilon_{1}^{\frac{1}{2}- \frac{\beta}{4} } E_{1}^{ \frac{\beta}{4} }.
\end{align}
Rewrite (4.3) as
\begin{align}
\partial_{u} \left[ r ( \partial_{v} \Psi )^{2} \right]
+ \partial_{v} \left[ r ( \partial_{u} \Psi )^{2} \right] = 2 G_{1}.
\end{align}
Take $ \alpha = \beta $ in \eqref{11} and rewrite it as

\begin{align}
&	-\partial_{v} \left( r^{1-\beta} q_0 q_1 \right)+\partial_v\left(\frac12 r^{1-\beta} (q_0^2+q_1^2)\right)-\frac 12 \partial_{v} (r^{-\beta+\sigma} Q_0 q_0)\nonumber\\
&
+\frac 12 \partial_v (q_1 r^{-\beta+\sigma} Q_0)+\frac{\beta-\sigma +1}{4} \partial_v (r^{2\sigma-1-\beta}   Q_0 ^2)\nonumber\\
&-\partial_{u} \left( r^{1-\beta} q_0 q_1 \right)-\partial_u\left(\frac12 r^{1-\beta} (q_0^2+q_1^2)\right)-\frac 12 \partial_{u} (r^{-\beta+\sigma} Q_0 q_0)\nonumber\\
&
-\frac 12 \partial_u (q_1 r^{-\beta+\sigma} Q_0)-\frac{\beta-\sigma +1}{4} \partial_u (r^{2\sigma-1-\beta}   Q_0 ^2)\nonumber\\
&	= \frac{q_0}{2r^{\beta-\sigma}}\int_0^r \xi^{-\sigma} A_0 q_1 d \xi- \frac{\sigma q_0}{2r^{\beta-\sigma}} \int _0^r \xi^{-\sigma-1} q_0(\xi) d \xi+ \frac 12  r^{-\beta+\sigma} A_0q_0  Q_0\nonumber\\
&-\frac{(1+\beta-2\sigma)(\beta-\sigma +1)}{4}r^{2\sigma-2-\alpha} Q_0^2-\frac\beta2 \frac{ q_0^2+q_1^2}{r^{\alpha}}\nonumber\\.
\end{align}

\begin{align}
	&\partial_{v}\left( r^{1-\beta}\frac{1}{2}\left(q_0-q_1\right)^2+\frac{1}{2}r^{-\beta+\sigma} Q_0\left(q_1-q_0\right)+\frac{\beta-\sigma +1}{4} r^{2\sigma-1-\beta}   Q_0 ^2\right)\nonumber\\
	&-	\partial_{u}\left( r^{1-\beta}\frac{1}{2}\left(q_0+q_1\right)^2+\frac{1}{2}r^{-\beta+\sigma} Q_0\left(q_1+q_0\right)+\frac{\beta-\sigma +1}{4} r^{2\sigma-1-\beta}   Q_0 ^2\right)\nonumber\\
	&= G_{\beta},
\end{align}
where
\begin{align}
G_{\beta}:=& \frac{q_0}{2r^{\beta-\sigma}}\int_0^r \xi^{-\sigma} A_0 q_1 d \xi- \frac{\sigma q_0}{2r^{\beta-\sigma}} \int _0^r \xi^{-\sigma-1} q_0(\xi) d \xi+ \frac 12  r^{-\beta+\sigma} A_0q_0  Q_0\nonumber\\
&-\frac{(1+\beta-2\sigma)(\beta-\sigma +1)}{4}r^{2\sigma-2-\alpha} Q_0^2-\frac\beta2 \frac{ q_0^2+q_1^2}{r^{\alpha}}.
\end{align}

We integrate above balace law  over space time, and take the $\sigma$ sufficiently small. By hardy inequality, we are arrive at the following estimates

\begin{align}
	&\frac{(1+\beta-2\sigma)(\beta-\sigma +1)}{4}r^{2\sigma-2-\alpha} Q_0^2+\frac\beta2 \frac{ q_0^2+q_1^2}{r^{\alpha}}\nonumber\\
	&\leq 2\left(\frac{(1+\beta-2\sigma)(\beta-\sigma +1)}{4}r^{2\sigma-2-\alpha} Q_0^2 +\frac\beta2 \frac{ q_0^2+q_1^2}{r^{\alpha}} +\frac{\sigma q_0}{2r^{\beta-\sigma}} \int _0^r \xi^{-\sigma-1} q_0(\xi) d \xi \right)\nonumber\\
	&\leq 2\left(\frac{\sigma q_0}{2r^{\beta-\sigma}} \int _0^r \xi^{-\sigma}A_0 q_0(\xi) d \xi+ \frac 12  r^{-\beta+\sigma} A_0q_0  Q_0\right) + C\epsilon^{1-\frac{\beta}{2}}_{1}E^{\frac{\beta}{2}}_1\nonumber\\
\end{align}

Noting that (2.10) and energy small, the first two terms on the right can be absorbed by the left hand side, therefore
\begin{align}
\int_{0}^{T}  \int_{0}^{\infty} |G_{\beta}|(t,r) drdt \lesssim \epsilon_{1}^{1- \frac{\beta}{2} } E_{1}^{ \frac{\beta}{2} }.
\end{align}

Noting that
\begin{align}
	&r^{1-\beta} \lvert  \Phi_{u}\rvert^2 =\frac{1}{4} r^{1-\beta}\left(q_0-q_1\right)^2\nonumber\\
	&\lesssim r^{1-\beta}\frac{1}{2}\left(q_0-q_1\right)^2+\frac{1}{2}r^{-\beta+\sigma} Q_0\left(q_1-q_0\right)+\frac{\beta-\sigma +1}{4} r^{2\sigma-1-\beta}   Q_0 ^2
\end{align}
and
\begin{align}
&r^{1-\beta} \lvert  \Phi_{v}\rvert^2 =\frac{1}{4} r^{1-\beta}\left(q_0-q_1\right)^2\nonumber\\
&\lesssim r^{1-\beta}\frac{1}{2}\left(q_0+q_1\right)^2+\frac{1}{2}r^{-\beta+\sigma} Q_0\left(q_1+q_0\right)+\frac{\beta-\sigma +1}{4} r^{2\sigma-1-\beta}   Q_0 ^2
\end{align}

We apply Lemma 3.2 to get
\begin{align}
\int_{0}^{T}  \int_{0}^{\infty} r^{2-\beta} \left[ ( \partial_{u} \Phi )^{2} ( \partial_{v} \Psi )^{2} + ( \partial_{u} \Psi )^{2} ( \partial_{v} \Phi )^{2} \right] drdt
\lesssim \epsilon_{1}^{1- \frac{\beta}{2} } E_{1}^{1+ \frac{\beta}{2} }.
\end{align}
So we get
\begin{align}
\int_{0}^{T}  \int_{0}^{\infty} |g_{1}| drdt
\lesssim \epsilon_{1} E_{1}.
\end{align}

Now we estimate $g_{2} $.
\begin{align}
& \int_{0}^{T}  \int_{0}^{\infty} |g_{2}| drdt
\lesssim \int_{0}^{T} \left( \int_{0}^{\infty} |\Psi_{t}|^{2} rdr \right)^{ \frac{1}{2} }
\left( \int_{0}^{\infty} |\Psi|^{2} |\Phi_{u}|^{2} |\Phi_{v}|^{2} rdr \right)^{ \frac{1}{2} } dt \notag \\
& \lesssim E_{1}^{\frac{1}{2} } \int_{0}^{T} \left( \int_{0}^{\infty} |\Psi|^{2} |\Phi_{u}|^{2} |\Phi_{v}|^{2} rdr \right)^{ \frac{1}{2} } dt.
\end{align}
We have
\begin{align*}
& \int_{0}^{\infty} |\Psi|^{2} |\Phi_{u}|^{2} |\Phi_{v}|^{2} rdr
= \int_{0}^{\infty} |\Phi_{u}|^{2} |\Phi_{v}|^{2} d\int_{0}^{r} |\Psi|^{2} \lambda d\lambda \\
& = -2\int_{0}^{\infty} \int_{0}^{r} |\Psi|^{2} \lambda d\lambda \left[  ( \Phi_{u} \cdot \hat{\Psi}_{u}  ) |\Phi_{v}|^{2}
+ |\Phi_{u}|^{2} ( \Phi_{v} \cdot \hat{\Psi}_{v}  ) \right]  \\
& \leq 2 \int_{0}^{\infty} dr \int_{0}^{r} |\Psi|^{2} \lambda^{- \beta} d\lambda \left[  r^{1+ \beta} |\Phi_{u}| |\Phi_{v}| ( |\hat{\Psi}_{u}| |\Phi_{v}|+ |\hat{\Psi}_{v}| |\Phi_{u}| ) \right]  \\
& \leq 2\left( \int_{0}^{\infty} |\Psi|^{2} r^{- \beta} dr \right) \left( \int_{0}^{\infty} (|\hat{\Psi}_{u}|^{2} |\Phi_{v}|^{2} +  |\hat{\Psi}_{v}|^{2} |\Phi_{u}|^{2}   ) r^{2- \beta} dr  \right)^{ \frac{1}{2} } \\
& \left( \int_{0}^{\infty} |\Phi_{u}|^{2} |\Phi_{v}|^{2}   r^{3 \beta} dr  \right)^{ \frac{1}{2} }.
\end{align*}
Therefore we get
\begin{align}
& \int_{0}^{T}  \int_{0}^{\infty} |g_{2}| drdt
\lesssim E_{1}^{ \frac{1}{2} }
 \left( \int_{0}^{T} \int_{0}^{\infty} |\Psi|^{2} r^{- \beta} drdt \right)^{ \frac{1}{2} } \notag \\
& \left( \int_{0}^{T} \int_{0}^{\infty} r^{2- \beta} (|\hat{\Psi}_{u}|^{2} |\Phi_{v}|^{2} +  |\hat{\Psi}_{v}|^{2} |\Phi_{u}|^{2}   ) drdt \right)^{ \frac{1}{4} } \notag \\
& \left( \int_{0}^{T} \int_{0}^{\infty} |\Phi_{u}|^{2} |\Phi_{v}|^{2}   r^{3 \beta} drdt  \right)^{ \frac{1}{4} }.
\end{align}
The first term in the right hand side can be estimated by (4.14), and the second term can be estimated in a way similar to (4.23). Thus, it remains to prove
\begin{align*}
\int_{0}^{T} \int_{0}^{\infty} |\Phi_{u}|^{2} |\Phi_{v}|^{2}   r^{3 \beta} drdt
\lesssim \epsilon_{1}^{1+\frac{3\beta}{2} } E_{1}^{ 1- \frac{3\beta}{2} }.
\end{align*}
We have by (4.17) with $\beta$ replaced by $\alpha$
\begin{align}
&\partial_{v}\left( r^{1-\alpha}\frac{1}{2}\left(q_0-q_1\right)^2+\frac{1}{2}r^{-\alpha+\sigma} Q_0\left(q_1-q_0\right)+\frac{\alpha-\sigma +1}{4} r^{2\sigma-1-\alpha}   Q_0 ^2\right)\nonumber\\
&-	\partial_{u}\left( r^{1-\alpha}\frac{1}{2}\left(q_0+q_1\right)^2+\frac{1}{2}r^{-\alpha+\sigma} Q_0\left(q_1+q_0\right)+\frac{\alpha-\sigma +1}{4} r^{2\sigma-1-\alpha}   Q_0 ^2\right)\nonumber\\
&= G_{\alpha},
\end{align}
We have
\begin{align*}
\Phi_{t} \cdot \left[   \Phi_{tt} - \Phi_{rr} - \frac{\Phi_{r}}{r} \right]  = 0,
\end{align*}
so we get
\begin{align*}
\partial_{u} \left[ r^{1-\alpha} ( \partial_{v} \Phi )^{2} \right]
+ \partial_{v} \left[ r^{1-\alpha} ( \partial_{u} \Phi )^{2} \right]  = \hat{G}_{\alpha},
\end{align*}
where
\begin{align*}
 \hat{G}_{\alpha} = (1-\alpha) r^{-\alpha} \left[ ( \partial_{v} \Phi )^{2} - ( \partial_{u} \Phi )^{2}  \right] .
\end{align*}
Thus, we have
\begin{align*}
\int_{0}^{T} \int_{0}^{\infty} \left( | \hat{G}_{\alpha} | +|G_{\alpha}|  \right) drdt
 \lesssim \epsilon_{1}^{1-\frac{ \alpha}{2} } E_{1}^{ \frac{ \alpha}{2} }.
\end{align*}
So we apply Lemma3.2 to get
\begin{align*}
\int_{0}^{T} \int_{0}^{\infty} r^{2-2\alpha} ( \partial_{v} \Phi )^{2} ( \partial_{u} \Phi )^{2}
 drdt
\lesssim \epsilon_{1}^{2-\alpha} E_{1}^{\alpha }.
\end{align*}
Take $ 2-2\alpha=3\beta $, we get
\begin{align*}
\int_{0}^{T} \int_{0}^{\infty} r^{3\beta} ( \partial_{v} \Phi )^{2} ( \partial_{u} \Phi )^{2}
drdt
\lesssim \epsilon_{1}^{1+\frac{3\beta}{2} } E_{1}^{ 1- \frac{3\beta}{2} }.
\end{align*}
This completes the proof of the regularity with small energy.

\begin{remark}
This paper was published in Chinese Annals of Mathematics in 2022, unfortunately there is a miss calculation in the main equality. The present version correct this mistakes.
\end{remark}

\bigskip

{\bf Acknowledgement} The author is supported by Key Laboratory of Mathematics for Nonlinear Sciences (Fudan University), Ministry of Education of China, P.R.China. Shanghai Key Laboratory for Contemporary Applied Mathematics, School of Mathematical Sciences, Fudan University, P.R. China, and by Shanghai Science and Technology Program [Project No. 21JC1400600].

\bibliographystyle{plain}
\end{CJK*}

\end{document}